\documentclass[12pt,a4paper,reqno]{amsart}
\usepackage{amssymb}
\usepackage{latexsym}
\usepackage{amsmath}
\usepackage{mathrsfs}
\usepackage{cite}

\usepackage{calc}

\newcommand{\scal}[2]{\langle #1,#2\rangle}
\newcommand{\rr}[1]{\mathbf R^{#1}}
\newcommand{\cc}[1]{\mathbf C^{#1}}
\newcommand{\nm}[2]{\Vert #1\Vert _{#2}}
\newcommand{\nmm}[1]{\Vert #1\Vert }

\newcommand{\op}{\operatorname{Op}}
\newcommand{\sets}[2]{\{ \, #1\, ;\, #2\, \} }

\newcommand{\fy}{\varphi}
\newcommand{\cdo}{\, \cdot \, }

\newcommand{\eabs}[1]{\langle #1\rangle}     

\newcommand{\hatconv}{\, \widehat *\, }
\newcommand{\tp}{\operatorname{Tp}}
\newcommand{\vrum}{\vspace{0.1cm}}
\newcommand{\im}{i}

\numberwithin{equation}{section}          

\newtheorem{thm}{Theorem}
\numberwithin{thm}{section}
\newtheorem{prop}[thm]{Proposition}
\newtheorem{cor}[thm]{Corollary}
\newtheorem{lemma}[thm]{Lemma}

\newcommand{\rubrik}{}

\theoremstyle{definition}

\newtheorem{defn}[thm]{Definition}

\theoremstyle{remark}
\newtheorem{rem}[thm]{Remark}

\newcommand{\rd}{\mathbf{R} ^{d}}

\newcommand{\intrd}{\int _{\rd }}

\title{Mapping properties for the Bargmann transform on modulation
spaces}

\author{Mikael Signal}

\address{Department of Mathematics, University of Agder, Kristiansand,
Norway}

\email{mikael.signal@uia.no}

\author{Joachim Toft}

\address{Department of Computer science, Physics and Mathematics,
Linn{\ae}us University, V{\"a}xj{\"o}, Sweden\footnote{Former address:
Department of Mathematics and Systems Engineering,
V{\"a}xj{\"o} University, Sweden}}

\email{joachim.toft@lnu.se}

\begin{document}

\begin{abstract}
We investigate mapping properties for the Bargmann transform and prove
that this transform is isometric and bijective from modulation spaces
to convenient Banach spaces of analytic functions.
\end{abstract}

\maketitle

\par

\section{Introduction}\label{sec0}

\par

In \cite {B1}, V. Bargmannn introduce a transform $\mathfrak V$ which
is bijective and isometric from $L^2(\rr d)$ into the Hilbert space
$A^2(\cc d)$ of all entire analytic functions $F$
on $\cc d$ such that $F\cdot e^{-|\cdo |^2/2}\in L^2(\cc d)$. (We use
the usual notations for the usual function
and distribution spaces. See e.{\,}g. \cite{Ho1}, and refer to Section
\ref{sec1} for specific definitions and other notations.) Furthermore,
several important
properties for $\mathfrak V$ were established. For example:
\begin{itemize}
\item the Hermite functions are mapped into the normalized analytical
monomials. Furthermore, the latter set forms an orthonormal basis for
$A^2(\cc d)$;

\vrum

\item the creation and annihilation operators, and harmonic oscillator
on appropriate elements in $L^2$, are transformed by $\mathfrak V$
into simple operators;

\vrum

\item there is a convenient reproducing formula for elements in $A^2$.
\end{itemize}

\par

In \cite{B2}, Bargmann continued his work and discussed mapping
properties for $\mathfrak V$ on more general spaces. For example,
he proves that $\mathfrak V(\mathscr S')$, the image of $\mathscr S'$
under the Bargmann transform is given by the formula
\begin{equation}\label{Sprimbild}
\mathfrak V(\mathscr S') = \cup _{\omega \in \mathscr
P}A^{\infty ,\infty}_{(\omega )},
\end{equation}
Here $A^{p,q}_{(\omega )}(\cc d)$ is the set of all entire functions
$F$ on $\cc d$ such that $F\cdot e^{-|\cdo
|^2/2}\cdot \omega _0$ belongs to the mixed Lebesgue space
$L^{p,q}(\cc d)$, for some appropriate modification
$\omega _0$ of the weight function $\omega$.

\medspace

The Bargmann transform can easily be reformulated in terms of the
short-time Fourier transform, with a particular Gauss function as
window function. In this context we remark that the (classical)
modulation spaces $M^{p,q}$, $p,q \in [1,\infty]$, as introduced by
Feichtinger in \cite{F1}, consist of all tempered distributions whose
short-time Fourier transforms (STFT) have finite mixed $L^{p,q}$
norm. It follows that the parameters $p$ and $q$ to some extent
quantify the degrees of asymptotic decay and singularity of
the distributions in $M^{p,q}$. The theory of modulation spaces was
developed further and generalized in
\cite{FG1,FG2,Feichtinger5,Gc1}, where
Feichtinger and Gr{\"o}chenig established the theory of coorbit
spaces. In particular, the modulation space $M^{p,q}_{(\omega )}$,
where $\omega$ denotes a weight function on
phase (or time-frequency shift) space, appears as the set of tempered
(ultra-) distributions
whose STFT belong to the weighted and mixed Lebesgue space
$L^{p,q}_{(\omega )}$. Here the weight $\omega$ quantifies the degree
of asymptotic decay and singularity of the distribution in
$M^{p,q}_{(\omega )}$.

\par

A major idea behind the design of these spaces was to find
useful Banach spaces, which are defined in a way similar to Besov
spaces, in the sense of replacing the dyadic decomposition on the
Fourier transform side, characteristic to Besov spaces, with a
\textit{uniform} decomposition. From the construction of these
spaces, it turns out that modulation spaces and Besov spaces in some
sense are rather similar, and sharp embeddings between these
spaces can be found in \cite{Toft2, To04B}, which
are improvements of certain embeddings in \cite {Grobner}. (See
also \cite {Sugimoto1} for verification of the sharpness.)

\par

During the last 15 years many results have been proved which confirm
the usefulness of the modulation
spaces in time-frequency analysis, where they occur naturally. For
example, in \cite{Feichtinger5,Gc2,GL}, it
is shown that all modulation spaces admit
reconstructible sequence space representations using Gabor frames.

\par

Parallel to this development, modulation spaces have been incorporated
into the calculus of pseudo-differential operators, which also involve
Toeplitz operators. (See
e.{\,}g. \cite{Gc2,GT,He1,HTW,Sugimoto1,Toft2,To04B,To5,To8} and the
references therein.)

\medspace

By reformulating the Bargmann transform in terms of the short-time Fourier
transform, and using the fundamental role of the short-time Fourier
transform in the definition of modulation spaces, it follows easily
that the Bargmann transform is continuous and injective from
$M^{p,q}_{(\omega )}$ to $A^{p,q}_{(\omega)}$. Furthermore, by
choosing the window function as a particular Gaussian function in the
$M^{p,q}_{(\omega )}$ norm, it follows that $\mathfrak V\, :\,
M^{p,q}_{(\omega)}\to A^{p,q}_{(\omega)}$ is isometric.

\par

These facts and several other mapping properties for the Bargmann
transform on modulation spaces
were established and proved by Feichtinger,
Gr{\"o}chenig and Walnut in \cite{FG1, FGW, Gc1, GW}. 
In fact, here they state and motivate that the Bargmann transform from
$M^{p,q}_{(\omega )}$ to
$A^{p,q}_{(\omega )}$ is not only injective, but in fact
\emph{bijective}. In their proof of the surjectivity, they recall from
\cite{B1} that the Bargmann transform is bijective from $\mathscr
S_0(\rr d)$ to $P(\cc d)$, where $\mathscr S_0(\rr d)$ is the set of
finite linear combinations of the Hermite functions and $P(\cc d)$ is
the set of analytic polynomials on $\cc d$.  Then they use duality in
combination with the argument that $P(\cc d)$ is dense in
$A^{p,q}_{(\omega )}$ when $p,q<\infty$. Since $\mathscr S_0$ is dense
in $M^{p,q}_{(\omega )}$, the asserted surjectivity easily follows
from these arguments.

\par

We are convinced that, somewhere in the litterature, it is proved that
$P(\cc d)$ is dense in $A^{p,q}_{(\omega )}$ (for example, a proof
might occur in \cite{FG1, Gc1, GW, JPR}). On the other hand, so far we
are unable to find any explicit proof of this fact. Especially, we
could not find any explicit references in the papers \cite{FG1, Gc1,
GW}.

\par

In \cite{FGW,Gc1,GW}, Feichtinger, Gr{\"o}chenig and Walnut also give
another motivation for the surjectivity. More precisely, they use the
arguments that the Bargmann-Fock repsresentation of the Heisenberg
group is unitarily equivalent to the Schr{\"o}dinger representation
with $\mathfrak V$ as the intertwining operator.
Then they explain that the general intertwining theorem \cite[Theorem
4.8]{FG1} applied to the Schr{\"o}dinger representation and the
Bargmann-Fock representation implies that $\mathfrak V$ extends to a
Banach space isomorphism from $M^{p,p}_{(\omega )}$ to
$A^{p,p}_{(\omega )}$, and the asserted surjectivity follows.

\par

It is obvious that these arguements are sufficent to conclude that
$\mathfrak V$ is a Banach space isomorphism from $M^{p,q}_{(\omega )}$
to $\mathfrak V(M^{p,q}_{(\omega )})$. On the other hand, so far we
are unable to understand that these arguments are sufficient to
conclude that indeed $\mathfrak V(M^{p,q}_{(\omega )})$ is the same
as $A^{p,q}_{(\omega )}$.

\medspace

In this paper we take an alternative approach for proving this
bijectivity. The main part is to prove that \eqref{Sprimbild} can be
improved into
\begin{equation}\tag*{(\ref {Sprimbild})$'$}
\mathfrak V(\mathscr S') = \cup _{\omega \in \mathscr
P}A^{p,q}_{(\omega )},
\end{equation}
when $p,q\in [1,\infty]$. Admitting this for a while, it follows that
each element in $A^{p,q}_{(\omega )}$ is a Bargmann transform of a
tempered distribution. The fact that the Bargmann transform is
continuous and injective from $M^{p,q}_{(\omega )}$ to
$A^{p,q}_{(\omega )}$ then shows that this tempered distribution must
belong to $M^{p,q}_{(\omega )}$, and the result follows.

\par

When proving \eqref{Sprimbild}$'$ we first consider mapping properties
on Hilbert spaces, defined by the  Harmonic oscillator. We prove that
such Hilbert spaces are modulation spaces of the form
$M^{2,2}_{(\omega )}$, when $\omega (x,\xi )=\sigma _N(x,\xi )=\eabs
{x,\xi }^N$ for some even number $N$. Here
$$
\eabs x=(1+|x|^2)^{1/2}\quad \text{and}\quad  \eabs {x,\xi
}=(1+|x|^2+|\xi |^2)^{1/2}
$$
as usual. Furthermore, we use the analysis
in \cite{B1,B2} to prove that $\mathfrak V$ maps $M^{2,2}_{(\sigma
_N)}$ bijectively and isometrically onto $A^{2,2}_{(\sigma _N)}$. Since
any tempered distribution belongs to $M^{2,2}_{(\sigma _N)}$ provided
$N$ is a large enough negative number, \eqref{Sprimbild}$'$ follows in
the case $p=q=2$.

\par

By using an argument of harmonic mean values, we thereafter prove that
similar facts also holds for $p=q=1$. Since $A^{p,q}_{(\omega
)}\subseteq A^{1,1}_{(\sigma _N)}$, by H{\"o}lder's inequality,
provided $N$ is a large enough negative number, it follows that each
element in $A^{p,q}_{(\omega )}$ is a Bargmann transform of a tempered
distribution. The asserted bijectivity is now a consequence of  the
fact that $\mathfrak V\, :\, M^{p,q}_{(\omega )}\to A^{p,q}_{(\omega
)}$ is continuous and injective.

\medspace

The paper is organized as follows. In Section \ref{sec1} we recall
some facts for modulation spaces and the Bargmann transform. In
Section \ref{sec2} we prove the main result, i.{\,}e. that the
Bargmann transform is bijective from $M^{p,q}_{(\omega )}$ to
$A^{p,q}_{(\omega )}$. In fact, we prove a more general result involving
general modulation spaces $M(\omega ,\mathscr B)$, parameterized with
the weight funciton $\omega$ and the translation invariant BF-space
$\mathscr B$. Finally, in Section \ref{sec3} we present some
consequences of the main result. Several of these consequences can be
found in  \cite{FG1, FGW, Gc1, GW}. However in our approach, such known
consequences enter the theory in different ways comparing to
\cite{FG1, FGW, Gc1, GW}.

\par

\section{Preliminaries}\label{sec1}

\par

In this section we give some definitions and recall some basic
facts. The proofs are in general omitted.

\par

\subsection{Translation invariant BF-spaces}\label{subsec1.1}

\par

We start with presenting appropriate conditions of the involved weight
functions. Assume that $\omega, v\in L^\infty _{loc}(\rr d)$ are
positive functions. Then $\omega$ is called $v$-moderate if
\begin{equation}\label{moderate}
\omega (x+y) \leq C\omega (x)v(y)
\end{equation}
for some constant $C$ which is independent of $x,y\in \rr d$. If $v$
in \eqref{moderate} can be chosen as a polynomial, then $\omega$ is
called polynomially moderate. We let $\mathscr P(\rr d)$ be the set
of all polynomially moderated functions on $\rr d$. We also let
$\mathscr P_0(\rr d)$ be the set of all $\omega \in \mathscr P(\rr
d)\cap C^\infty (\rr d)$ such that $(\partial ^\alpha \omega )/\omega
\in L^\infty$, for every multi-index $\alpha$. We remark that for each
$\omega \in \mathscr P(\rr d)$, there is an element $\omega _0\in
\mathscr P_0(\rr d)$ which is equivalent to $\omega$, in the sense
that for some constant $C$, it holds
$$
C^{-1}\omega _0\le \omega \le C\omega _0
$$
(cf. e.{\,}g. \cite{To04B,To5}).

\par

We say that $v$ is
\emph{submultiplicative} when \eqref{moderate} holds with $\omega =v$.
Throughout we assume that the submultiplicative weights are
even. Furthermore, $v$ and $v_j$ for $j\ge 0$, always stand for
submultiplicative weights if nothing else is stated.

\par

An important type of weight functions is 
\begin{equation}\label{defsigmas}
\sigma _s(x)\equiv \eabs {x}^s = (1+|x|^2)^{s/2},
\end{equation}

\par

For each $\omega \in \mathscr P(\rr d)$ and $p\in [1,\infty ]$, we let
$L^p_{(\omega )}(\rr d)$ be the Banach space which consists of all
$f\in L^1_{loc}(\rr d)$ such that $\nm {f}{L^p_{(\omega )}}\equiv \nm
{f\, \omega }{L^p}$ is finite.

\par

Next we recall the definition of Banach function spaces (BF-spaces).

\par

\begin{defn}\label{BFspaces}
Assume that $\mathscr B$ is a Banach space of complex-valued
measurable functions on $\rr d$ and that $v \in \mathscr P(\rr {d})$
is submultiplicative.
Then $\mathscr B$ is called a \emph{(translation) invariant
BF-space on $\rr d$} (with respect to $v$), if there is a constant $C$
such that the following conditions are fulfilled:
\begin{enumerate}
\item $\mathscr S(\rr d)\subseteq \mathscr
B\subseteq \mathscr S'(\rr d)$ (continuous embeddings).

\vrum

\item If $x\in \rr d$ and $f\in \mathscr B$, then $f(\cdot -x)\in
\mathscr B$, and
\begin{equation}\label{translmultprop1}
\nm {f(\cdot -x)}{\mathscr B}\le Cv(x)\nm {f}{\mathscr B}\text .
\end{equation}

\vrum

\item if  $f,g\in L^1_{loc}(\rr d)$ satisfy $g\in \mathscr B$  and $|f|
\le |g|$ almost everywhere, then $f\in \mathscr B$ and
$$
\nm f{\mathscr B}\le C\nm g{\mathscr B}\text .
$$

\vrum

\item the map $(f,\fy )\mapsto f*\fy$ is continuous from $\mathscr
B\times C_0^\infty (\rr d)$ to $\mathscr B$.
\end{enumerate}
\end{defn}

\par

\begin{rem}\label{0303}
Assume that $\mathscr B$ is a translation invariant BF-space. If $f\in
\mathscr B$ and $h\in L^\infty$, then it follows from (3) in
Definition \ref{BFspaces} that $f\cdot h\in \mathscr B$ and
\begin{equation}\label{multprop}
\nm {f\cdot h}{\mathscr B}\le C\nm f{\mathscr B}\nm h{L^\infty}.
\end{equation}
\end{rem}

\par

\begin{rem}\label{newbfspaces}
Assume $\omega _0,v,v_0\in \mathscr P(\rr d)$ are such that $v$ and
$v_0$ are submultiplicative,
$\omega _0$ is $v_0$-moderate, and assume that $\mathscr B$ is a
translation-invariant BF-space on $\rr d$ with respect to $v$. Also
let $\mathscr B(\omega _0)$ be the Banach space which consists of all
$f\in L^1_{loc}(\rr d)$ such that $\nm f{\mathscr B(\omega _0)}\equiv
\nm {f\, \omega _0}{\mathscr B}$ is finite. Then $\mathscr B(\omega
_0)$ is a translation invariant BF-space with respect to $v_0v$.
\end{rem}

\par

\begin{rem}\label{convext}
Let $\mathscr B$ be a translation invariant BF-space on $\rr d$ with
respect to $v\in \mathscr P(\rr d)$. Then
it follows that the map $(f,g )\mapsto f*g$ from $\mathscr
B\times C_0^\infty (\rr d)$ to $\mathscr B$ extends uniquely to a
continuous mapping from $\mathscr B\times L^1_{(v)} (\rr d)$ to
$\mathscr B$. In fact, if $f\in \mathscr B$ and $g\in C_0^\infty (\rr d)$, then Minkowski's inequality gives
\begin{multline*}
\nm {f*g}{\mathscr B} =\Big \Vert \int f(\cdo -y)g(y)\, dy \Big \Vert
_{\mathscr B}
\\[1ex]
\le \int \nm { f(\cdo-y)}{\mathscr B}|g(y)|\, dy \le C\int \nm {
f}{\mathscr B}|g(y)v(y)|\, dy = C\nm { f}{\mathscr B}\nm
g{L^1_{(v)}}.
\end{multline*}
The assertion is now a consequence of the fact that $C_0^\infty$ is
dense in $L^1_{(v)}$.
\end{rem}

\par

\begin{rem}\label{BFemb}
Let $\mathscr B$ be an invariant BF-space. Then it is easy to find
Sobolev type spaces which are continuously embedded in $\mathscr
B$. In fact, for each $p\in [1,\infty )$ and integer $N\ge 0$, let
$Q^p_N(\rr d)$ be the set of all $f\in L^p(\rr d)$ such that $\nm
f{Q^p_N}<\infty$, where
$$
\nm f{Q^p_N}\equiv \sum _{|\alpha +\beta |\le N}\nm {x^\alpha D^\beta f}{L^p}.
$$
Then for each $p$ fixed, the topology for $\mathscr S(\rr d)$ can be
defined by the semi-norms $f\mapsto \nm f{Q^p_N}$, for
$N=0,1,\dots$. Since $\mathscr S$ is continuously embedded in the
Banach space $\mathscr B$, it now follows that
$$
\nm f{\mathscr B}\le C_{N,p}\nm f{Q^p_N}
$$
for some constants $C$ and $N$ which are independent of $f\in \mathscr
S$. Consequently, if in addition $p<\infty$, then $Q^p_N(\rr
d)\subseteq \mathscr B$, since $\mathscr S$ is dense in $Q^p_N$. This
proves the assertion.
\end{rem}

\par

The following proposition shows that even stronger embeddings
comparing to $Q^p_N(\rr d)\subseteq \mathscr B$ in
Remark \ref{BFemb} hold when $p=\infty$. Here, we set $L^p_N=L^p_{(\omega)}$
when $\omega (x)=\langle x\rangle^N$.

\par

\begin{prop}\label{Binvemb}
Let $\mathscr B$ be a translation invariant BF-space on $\rr d$, and
let $\omega \in \mathscr P(\rr d)$. Then
there is a large number $N$ such that
$$
L^\infty _N(\rr d)\subseteq \mathscr B(\omega )\subseteq L^1_{-N}(\rr
d).
$$
\end{prop}

\par

\begin{proof}
We may assume that $\omega =1$ in view of Remark \ref{newbfspaces}.
First we note that $\nm {\eabs \cdo ^{-N}}{\mathscr B}<\infty$,
provided $N$ is large enough. In fact, since $\mathscr S$ is
continuously embedded in $\mathscr B$, it follows that
\begin{equation}\label{BSest}
\nm f{\mathscr B}\le C\sum _{|\beta |\le N}\nm {\eabs \cdo ^N(\partial
^\beta f)}{L^\infty},
\end{equation}
when $f\in \mathscr S$, for some choices of constants $C$ and $N$, and
the assertion now follows since the right-hand side of \eqref{BSest}
is finite when $f(x)=\eabs x^{-N}$.

\par

It follows from Definition \ref{BFspaces} (2) that
\begin{equation}\label{translest}
\eabs x^{-d-1}\nm f{\mathscr B}\ge C\eabs x^{-N}\nm {f(\cdo
-x)}{\mathscr B},
\end{equation}
for some constants $C$ and $N$.  
By integrating \eqref{translest} and using Minkowski's inequality we
get
\begin{multline*}
\nm f{\mathscr B} \ge C_1\int \eabs x^{-N}\nm {f(\cdo
-x)}{\mathscr B}\, dx
\\[1ex]
\ge
C_2 \Big \Vert \int \eabs x^{-N}|f(\cdo -x)|\, dx\Big \Vert _{\mathscr
B}
=
C_2 \Big \Vert \int \eabs {x-\cdo }^{-N}|f(x)|\, dx\Big \Vert
_{\mathscr B}
\\[1ex]
\ge
C_3 \Big \Vert \int \eabs {x}^{-N}|f(x)|\, dx\, \eabs{\cdo }^{-N}\Big
\Vert _{\mathscr B}=C \nm f{L^1_{-N}},
\end{multline*}
for some positive constants $C_1,\dots ,C_3$, where $C=C_3\nm {\eabs
\cdo ^{-N}}{\mathscr B}<\infty$. This proves $\mathscr B\subseteq
L^1_{-N}$.

\par

It remains to prove $L^\infty _N\subseteq \mathscr B$. Let $N$ be as
in the first part of the proof. By straight-forward computations and
using remark \ref{0303} we
get
$$
\nm f{\mathscr B} = \nm {(f\eabs \cdo ^N)\eabs \cdo ^{-N}}{\mathscr
B}\le C_1\nm {\eabs \cdo ^{-N}}{\mathscr B}\nm f{L^\infty _N} =C\nm
f{L^\infty _N},
$$
for some constant $C_1$, where $C=C_1\nm {\eabs \cdo ^{-N}}{\mathscr
B}<\infty$. Hence $L^\infty _N\subseteq \mathscr B$, and the
result follows.
\end{proof}

\par

\subsection{The short-time Fourier transform and Toeplitz operators}

\par

Before giving the definition of short-time Fourier transform we recall some properties for the (usual) Fourier transform. The Fourier transform $\mathscr F$ is the linear and continuous
mapping on $\mathscr S'(\rr d)$ which takes the form
$$
(\mathscr Ff)(\xi )= \widehat f(\xi ) \equiv (2\pi )^{-d/2}\int _{\rr
{d}} f(x)e^{-i\scal  x\xi }\, dx
$$
when $f\in L^1(\rr d)$. Here $\scal \cdo \cdo$ denotes the usual scalar product on $\rr d$. The map $\mathscr F$ is a homeomorphism
on $\mathscr S'(\rr d)$ which restricts to a homeomorphism on
$\mathscr S(\rr d)$ and to a unitary operator on $L^2(\rr d)$.

\par

Let $\phi \in \mathscr S(\rr d)\setminus 0$ be fixed. For every $f\in
\mathscr S'(\rr d)$, the \emph{short-time Fourier transform} $V_\phi
f$ is the distribution on $\rr {2d}$ defined by the formula
\begin{equation}\label{defstft}
(V_\phi f)(x,\xi ) =\mathscr F(f\, \overline{\phi (\cdo -x)})(\xi ).
\end{equation}
We note that the right-hand side defines an element in $\mathscr
S'(\rr {2d})\cap C^\infty (\rr {2d})$.
We also note that if $f\in L^q_{(\omega )}$ for some $\omega \in
\mathscr P(\rr d)$, then $V_\phi f$ takes the form
\begin{equation}\tag*{(\ref{defstft})$'$}
V_\phi f(x,\xi ) =(2\pi )^{-d/2}\int _{\rr d}f(y)\overline {\phi
(y-x)}e^{-i\scal y\xi}\, dy.
\end{equation}

\medspace

Next we consider Toeplitz operators, also known as localization operators.
If $a\in \mathscr S'(\rr {2d})$ and $\phi \in \mathscr S(\rr d)\setminus 0$ are fixed, then the \emph{Toeplitz operator} $\tp (a)=\tp _\fy (a)$ is the linear and continuous operator on $\mathscr S(\rr d)$, defined by the formula
\begin{equation}\label{toeplitzdef}
(\tp _\phi (a)f,g)_{L^2(\rr d)} = (a\, V_\phi f,V_\phi g)_{L^2(\rr {2d})}.
\end{equation}
There are several characterizations of Toeplits operators and several ways to extend the definition of such operators (see e.{\,}g. \cite {GT} and the references therein). For example, the definition of $\tp _\phi (a)$ is uniquely extendable to every $a\in \mathscr S'(\rr {2d})$, and then $\tp _\phi (a)$ is still continuous on $\mathscr S(\rr d)$ to $\mathscr S(\rr d)$, and uniquely extendable to a continuous operator on $\mathscr S'(\rr d)$.

\par

Toeplitz operators  arise in pseudo-differential calculus~\cite{Fo,Le}, in the theory of quantization (Berezin
quantization~\cite{Berezin71}), and in signal
processing~\cite{Dau}   (under the name of 
time-frequency localization operators or STFT multipliers).

\par

\subsection{Modulation spaces}

\par

We shall now discuss modulation spaces and recall some basic properties. We start with the following definition.

\par

\begin{defn}\label{bfspaces2}
Let $\mathscr B$ be a translation invariant BF-space on
$\rr {2d}$, $\omega \in \mathscr P(\rr {2d})$, and let $\phi \in
\mathscr S(\rr d)\setminus 0$. Then the \emph{modulation space}
$M(\omega ,\mathscr B)$ consists of all $f\in
\mathscr S'(\rr d)$ such that
\begin{equation}\label{modnorm2}
\nm f{M(\omega ,\mathscr B)}
\equiv \nm {V_\phi f\, \omega }{\mathscr B}<\infty .
\end{equation}
If $\omega =1$, then the notation $M(\mathscr B)$ is used
instead of $M(\omega ,\mathscr B)$.
\end{defn}

\par

We remark that the modulation space $M(\omega ,\mathscr B)$ in
Definition \ref{bfspaces2} is independent of the choice of $\phi \in
\mathscr S(\rr d)\setminus 0$, and different choices of $\phi$ give
rise to equivalent norms.

\par

An important case concerns when $\mathscr B$ is a mixed-norm space of Lebesgue type. More precisely, let $\omega \in \mathscr P(\rr {2d})$, $p,q\in
[1,\infty ]$, and let $L^{p,q}(\rr {2d})$ be the
Banach space which consists of all  $F\in
L^1_{loc}(\rr {2d})$ such that
$$
\nm F{L^{p,q}} \equiv \Big (\intrd \Big (\intrd |F(x,\xi
 )|^p\, dx\Big )^{q/p}\, d\xi \Big )^{1/q}<\infty \, .
$$
(with obvious modifications when $p=\infty$ or
$q=\infty$). Also let $L^{p,q}_{*}(\rr {2d})$ be the set of all  $F\in
L^1_{loc}(\rr {2d})$ such that
$$
\nm F{L^{p,q}_{*}} \equiv \Big (\intrd \Big (\intrd |F(x,\xi
 )|^q\, d\xi \Big )^{p/q}\, dx \Big )^{1/p}<\infty \, .
$$
Then $M(\omega ,L^{p,q}(\rr {2d}))$ is the usual modulation space
$M^{p,q}_{(\omega )}(\rr d)$, and $M(\omega ,L^{p,q}_{*}(\rr {2d}))$ is
the space $W^{p,q}_{(\omega )}(\rr d)$ which is related to certain
types of classical Wiener amalgam spaces. For convenience we use the
notation $M^p_{(\omega)}$ instead of $M^{p,p}_{(\omega
)}=W^{p,p}_{(\omega )}$.

\par

For conveniency we set $M^{p,q}_s=M^{p,q}_{(\sigma _s)}$ and
$M^{p}_s=M^{p}_{(\sigma _s)}$, where $\sigma _s$ is given by \eqref{defsigmas}. Furthermore,  for  $\omega =1$ we set
$$
M(\mathscr B) = M(\omega ,\mathscr B),\quad M^{p,q} =
M^{p,q}_{(\omega )},\quad W^{p,q} =
W^{p,q}_{(\omega )},\quad \text{and}\quad  M^{p} = M^{p}_{(\omega )}
\, .
$$
Here we note that $\sigma _s$ depends on both $x$ and $\xi$-variables, which implies that
$$
\sigma _s(x,\xi )=\eabs {x,\xi }^s = (1+|x|^2+|\xi |^2)^{s/2}.
$$

\par

In the following proposition we recall some facts about
modulation spaces. We omit the proof, since the result can be found
in \cite{F1,FG1,FG2,Gc2,To5}.

\par

\begin{prop}\label{p1.4}
Let $p,q,p_j,q_j\in [1,\infty ]$, $\omega
,\omega _j,v,v_0\in \mathscr P(\rr {2d})$ for $j=1,2$ be such that
$\omega$ is $v$-moderate, and let $\mathscr B$ be a translation
invariant BF-space on $\rr {2d}$ with respect to $v_0$. Then the
following is true:
\begin{enumerate}
\item[{\rm{(1)}}] if $\phi \in M^1_{(v_0v)}(\rr d)\setminus 0$, then
$f\in M{(\omega ,\mathscr B)}$ if and only if \eqref{modnorm2} holds,
i.{\,}e. $M{(\omega ,\mathscr B)}$ is independent of the choice of
$\phi$. Moreover, $M{(\omega ,\mathscr B)}$ is a Banach space under the
norm in \eqref{modnorm2}, and different choices of $\phi$ give rise to
equivalent norms;

\vrum

\item[{\rm{(2)}}] if  $p_1\le p_2$,
$q_1\le q_2$ and $\omega _2\le C \omega _1$ for some constant $C$, then
\begin{alignat*}{2}
\mathscr S(\rr d) \subseteq M^{p_1,q_1}_{(\omega _1)}(\rr d)  
   &\subseteq M^{p_2,q_2}_{(\omega _2)}(\rr d)\subseteq
\mathscr S'(\rr d),
\\[1ex]
M^1_{(v_0v)}(\rr d)\subseteq M(&\omega ,\mathscr B)\subseteq
M^{\infty}_{(1/(v_0v))}(\rr d)\text ;
\end{alignat*}

\vrum

\item[{\rm{(3)}}] the sesqui-linear form $( \cdo ,\cdo )_{L^2}$ on
$\mathscr S(\rr d)$ extends to a continuous map from $M^{p,q}_{(\omega
)}(\rr d)\times M^{p'\! ,q'}_{(1/\omega )}(\rr d)$ to $\mathbf
C$. This extension is unique, except when $p=q'\in \{ 1,\infty \}$. On the
other hand, if $\nmm a = \sup |{(a,b)_{L^2}}|$, where the supremum is
taken over all $b\in M^{p',q'}_{(1/\omega )}(\rr d)$ such that
$\nm b{M^{p',q'}_{(1/\omega )}}\le 1$, then $\nmm {\cdot}$ and $\nm
\cdot {M^{p,q}_{(\omega )}}$ are equivalent norms;

\vrum

\item[{\rm{(4)}}] if $p,q<\infty$, then $\mathscr S(\rr d)$ is dense
in $M^{p,q}_{(\omega )}(\rr d)$, and the dual space of $M^{p,q}_{(\omega
)}(\rr d)$ can be identified with $M^{p'\! ,q'}_{(1/\omega )}(\rr
d)$, through the form $(\cdo  ,\cdo )_{L^2}$. Moreover,
$\mathscr S(\rr d)$ is weakly dense in $M^{\infty }_{(\omega )}(\rr
d)$.
\end{enumerate}
\end{prop}

\par

The following proposition is now a consequence of Remark 1.3 (5) in
\cite{To8} and Proposition \ref{p1.4} (2).

\par

\begin{prop}\label{capcupmodsp}
Let $\mathscr B$ be a translation invariant BF-space on $\rr {2d}$ and
let $\omega _j$ for $j\in J$ be a family of elements in $\mathscr
P(\rr {2d})$ such that for each $s\ge 0$, there is a constant
$C>0$, and $j_1,j_2\in J$ such that
$$
\omega _{j_1}(x,\xi )\le C\eabs {x,\xi}^{-s}\quad \text{and}\quad
C^{-1}\eabs {x,\xi}^{s} \le \omega _{j_2}(x,\xi ).
$$
Then
$$
\cup _{j\in J}M(\omega _{j},\mathscr B)=\mathscr S'(\rr d)\quad
\text{and}\quad \cap _{j\in J}M(\omega _{j},\mathscr B)=\mathscr S(\rr
d).
$$
\end{prop}

\par

\subsection{The Bargmann transform}

\par

We shall now consider the Bargmann transform which is defined by the
formula
\begin{equation}\label{bargtransf}
(\mathfrak Vf)(z) =\pi ^{-d/4}\int _{\rr d}\exp \Big ( -\frac 12(\scal
z z+|y|^2)+2^{1/2}\scal zy \Big )f(y)\, dy,
\end{equation}
when $f\in L^2(\rr d)$. We note that if $f\in
L^2(\rr d)$, then the Bargmann transform
$\mathfrak Vf$ of $f$ is the entire function on $\cc d$, given by
$$
(\mathfrak Vf)(z) =\int \mathfrak A_d(z,y)f(y)\, dy,
$$
or
\begin{equation}\label{bargdistrform}
(\mathfrak Vf)(z) =\scal f{\mathfrak A_d(z,\cdo )},
\end{equation}
where the Bargmann kernel $\mathfrak A_d$ is given by
$$
\mathfrak A_d(z,y)=\pi ^{-d/4} \exp \Big ( -\frac 12(\scal
zz+|y|^2)+2^{1/2}\scal zy\Big ).
$$
Here
$$
\scal zw = \sum _{j=1}^dz_jw_j,\quad \text{when} \quad z=(z_1,\dots ,z_d) \in \cc d\quad  \text{and} \quad w=(w_1,\dots ,w_d)\in \cc d,
$$
and $\scal \cdo \cdo $ denotes the duality between elements in $\mathscr S(\rr d)$ and $\mathscr S'(\rr d)$.
We note that the right-hand side in \eqref{bargdistrform} makes sense
when $f\in \mathscr S'(\rr d)$ and defines an element in $A(\cc d)$,
since $y\mapsto \mathfrak A_d(z,y)$ can be interpreted as an element
in $\mathscr S(\rr d)$ with values in $A(\cc d)$. Here and in what follows,
$A(\cc d)$ denotes the set of all entire functions on $\cc d$.

\par

From now on we assume that $\phi$ in \eqref{defstft},
\eqref{defstft}$'$  and \eqref{modnorm2} is given by
\begin{equation}\label{phidef}
\phi (x)=\pi ^{-d/4}e^{-|x|^2/2},
\end{equation}
if nothing else is stated. Then it follows that the Bargmann transform
can be expressed in terms of the short-time Fourier transform
$f\mapsto V_\phi f$. More precisely, for such choice of $\phi$, it
follows by straight-forward computations that
\begin{multline}\label{bargstft1}
  (\mathfrak{V} f)(z)  =  (\mathfrak{V}f)(x+\im \xi) 
 =  e^{(|x|^2+|\xi|^2)/2}e^{-i\scal x\xi}V_\phi f(2^{1/2}x,-2^{1/2}\xi
)
\\[1ex]
=e^{(|x|^2+|\xi|^2)/2}e^{-i\scal x\xi}(S^{-1}(V_\phi f))(x,\xi ),
\end{multline}
or equivalently,
\begin{multline}\label{bargstft2}
  V_\phi f(x,\xi )  =  
         e^{-(|x|^2+|\xi |^2)/4}e^{-\im \scal x \xi /2}(\mathfrak{V}f)
         (2^{-1/2}x,-2^{-1/2}\xi).
\\[1ex]
=e^{-i\scal x\xi /2}S(e^{-|\cdo |^2/2}(\mathfrak{V}f))(x,\xi ).
\end{multline}
Here $S$ is the dilation operator given by
\begin{equation}\label{Sdef}
(SF)(x,\xi ) = F(2^{-1/2}x,-2^{-1/2}\xi ).
\end{equation}
For future references we observe that \eqref{bargstft1} and \eqref{bargstft2} can be formulated into
$$
\mathfrak V = U_{\mathfrak V}\circ V_\phi ,\quad \text{and}\quad U_{\mathfrak V}^{-1} \circ \mathfrak V =  V_\phi ,
$$
where $U_{\mathfrak V}$ is the linear, continuous and bijective operator on $\mathscr D'(\rr {2d})=\mathscr D'(\cc d)$, given by
\begin{equation}\label{UVdef}
(U_{\mathfrak V}F)(x,\xi ) = e^{(|x|^2+|\xi |^2)/2}e^{-i\scal x\xi}F(2^{1/2}x,-2^{1/2}\xi ) .
\end{equation}

\par

We are now prepared to make the following definition.

\par

\begin{defn}\label{thespaces}
Let $\omega \in \mathscr P(\rr {2d})$ and let $\mathscr B$ be a
translation invariant BF-space on $\rr {2d}=\cc d$.
\begin{enumerate}
\item The space $\mathscr B_{\mathfrak V}(\omega )$ is the modified
weighted
$\mathscr B$-space which consists of all $F\in L^1_{loc}(\rr {2d})=
L^1_{loc}(\cc {d})$ such that
$$
\nm F{\mathscr B_{\mathfrak V}(\omega )}\equiv \nm {(S(Fe^{-|\cdo
|^2/2}))\omega}{\mathscr B}<\infty .
$$
Here $S$ is the dilation operator given by \eqref{Sdef};

\vrum

\item The space, $A(\omega ,\mathscr B)$ consists of all $F\in A(\cc
d)\cap \mathscr B_{\mathfrak V}(\omega )$ with topology inherited
from $\mathscr B_{\mathfrak V}(\omega )$;

\vrum

\item The space $A_0(\omega ,\mathscr B)$ is given by
$$
A_0(\omega ,\mathscr B) \equiv \sets {(\mathfrak Vf)}{f\in M(\omega
,\mathscr B)},
$$
and is equipped with the norm $\nm F{A_0(\omega ,\mathscr B)}\equiv \nm
f{M(\omega ,\mathscr B)}$, when $F=\mathfrak Vf$.
\end{enumerate}
\end{defn}

\par

The following result shows that the norm in $A_0(\omega,\mathcal{B})$ is well-defined.

\par

\begin{prop}\label{mainstep1prop}
Let $\omega \in \mathscr P(\rr {2d})$, let $\mathscr B$ be an
invariant BF-space on $\rr {2d}$ and let $\phi$ be as in
\eqref{phidef}. Then $A_0(\omega ,\mathscr
B)\subseteq A(\omega ,\mathscr B)$, and the map $\mathfrak V$ is an
isometric injection from $M(\omega ,\mathscr B)$ to $A(\omega ,\mathscr
B)$.
\end{prop}

\par

\begin{proof}
The result is an immediate consequence of \eqref{bargstft1},
\eqref{bargstft2} and Definition \ref{thespaces}.
\end{proof}

\par

We employ the same notational conventions for the spaces of type $A$ and $A_0$
 as we do for the modulation spaces.
In the case $\omega =1$ and $\mathscr B=L^2$, it follows from
\cite {B1} that Proposition \ref{mainstep1prop}
holds, and the inclusion is replaced by equality. That is, we have
$A^2_0=A^2$ which is called the Bargmann-Foch space, or just the Foch
space. In the next section we improve the latter property and show
that for any choice of $\omega \in \mathscr P$ and every translation
invariant BF-space $\mathscr B$, we have $A_0(\omega ,\mathscr B)=A(\omega
,\mathscr B)$.

\par

\section{Mapping results for the Bargmann
transform on modulation spaces}\label{sec2}

\par

In this section we prove that $A_0(\omega ,\mathscr B)$ is equal to
$A(\omega ,\mathscr B)$ for every choice of $\omega$ and $\mathscr
B$. That is, we have the following.

\par

\begin{thm}\label{mainthm}
Let $\mathscr B$ be a translation invariant BF-space on $\rr {2d}$ and
let $\omega \in \mathscr P(\rr {2d})$. Then $A_0(\omega ,\mathscr
B)=A(\omega ,\mathscr B)$, and the map $f\mapsto \mathfrak Vf$ from
$M(\omega ,\mathscr B)$ to $A(\omega ,\mathscr B)$ is isometric and
bijective.
\end{thm}

\par

We need some preparations for the proof, and start with giving some
remarks on the images of $\mathscr S(\rr d)$ and $\mathscr S'(\rr d)$
under the Bargmann transform. We denote these images by $A_{\mathscr
S}(\cc d)$ and $A_{\mathscr S}'(\cc d)$ respectively, i.{\,}e.
$$
A_{\mathscr S}(\cc d) \equiv \sets {\mathfrak Vf}{f\in \mathscr S(\rr
d)}\quad \text{and}\quad
A_{\mathscr S}'(\cc d) \equiv \sets {\mathfrak Vf}{f\in \mathscr
S'(\rr d)}.
$$
As a consequence of \eqref{bargstft1} and Propositions \ref{capcupmodsp} and \ref{mainstep1prop}, the inclusion
\begin{equation}\label{mapofSprime}
A_{\mathscr S}'(\cc d) \subseteq \sets {F\in A(\cc d)}{\nm {Fe^{-|\cdo
|^2/2}\sigma _{-N}}{L^p}<\infty\ \text{for some}\ N\ge 0}
\end{equation}
holds. We recall that in \cite{B2} it is proved that \eqref{mapofSprime}
holds with equality when $p=\infty$. An essential part of our
investigations concerns to prove that equality is attained in
\eqref{mapofSprime} for each $p\in [1,\infty ]$.

\par

\subsection{The image of harmonic oscillator on
$M^2_{2N}$.}\label{subsec2.1}

\par

Next we discuss mapping properties for a modified harmonic oscillator
on modulation spaces of the form $M^2_{2N}(\rr d)$, when $N$ is an
integer. The operator we have in mind is given by
\begin{equation}\label{harmosc}
H\equiv |x|^2 -\Delta +d+1,
\end{equation}
and we show that they are bijective between appropriate modulation spaces. Since Hermite functions constitute an orthonormal basis for $L^2=M^2$ and are eigenfunctions to the harmonic oscillator, we shall combine these facts to prove that dilations of such functions constitute an orthonormal basis for $M^2_{2N}$, for every integer $N$.

\par

We recall that if $\phi$ is given by \eqref{phidef} and
$$
a(x,\xi ) = \sigma _2(x,\xi )=|x|^2+|\xi |^2+1,
$$
then $H=\tp _\phi (a)$ (cf. e.{\,}g. Section 3 in \cite{To6}).
Let $\mathscr B$ be a translation invariant BF-space and $\omega \in
\mathscr P(\rr {2d})$. By Theorem 3.1 in \cite{GT} it now
follows that $H=\tp (a)=\tp (\sigma _2)$ is a continuous isomorphism
 from $M(\sigma _2\omega ,\mathscr B)$ to $M(\omega ,\mathscr B)$.
Since this holds for any weight $\omega$, induction
together with Banach's theorem now show that the following is true.

\par

\begin{prop}\label{identprop}
Let $N$ be an integer, $\omega \in \mathscr P(\rr {2d})$ and let
$\mathscr B$ be an invariant BF-space. Then $H^N$ on $\mathscr S'(\rr
d)$ restricts to a continuous isomorphism from $M(\sigma _{2N}\omega
,\mathscr B)$ to $M(\omega ,\mathscr B)$. In particular, the set
\begin{equation*}
\sets {f\in \mathscr S'(\rr d)}{H^Nf\in  L^2(\rr d)}
\end{equation*}
is equal to $M^2_{2N}(\rr d)$, and the norm $f\mapsto \nm {H^Nf}{L^2}$
is equivalent to $\nm f{M^2_{2N}}$.
\end{prop}

\par

From now on we assume that the norm
and scalar product of $M^2_{2N}(\rr d)$ are given by
$$
\nm f{M^2_{2N}}\equiv \nm {H^Nf}{L^2}\quad \text{and}\quad
(f,g)_{M^2_{2N}}\equiv (H^Nf,H^Ng)_{L^2}
$$
respectively. Then it follows from Proposition \ref{identprop} that $(e_j)_{j\in J}$ is an 
orthonormal basis for $M^2_{2N}$ if and only if $(H^Ne_j)_{j\in J}$ is an orthonormal basis for $L^2$. 
In the following we use this fact to find an
appropriate orthonormal basis for $M^2_{2N}(\rr d)$ in terms of 
Hermite functions.

\par

More precisely, we recall that the Hermite function $h_\alpha$ with respect to the multi-index $\alpha \in \mathbf N^d$ is defined by
$$
h_\alpha (x) = \pi ^{-d/4}(-1)^{|\alpha |}(2^{|\alpha |}\alpha
!)^{-1/2}e^{|x|^2/2}(\partial ^\alpha e^{-|x|^2}).
$$
The set $(h_\alpha )_{\alpha \in \mathbf N^d}$ is an orthonormal basis for $L^2$, and it follows from the definitions that  $h_\alpha $ is an
 eigenvector of $H$ with eigenvalue $2|\alpha
|+2d+1$ for every $\alpha$, i.{\,}e.  $Hh_\alpha = (2|\alpha
|+2d+1)h_\alpha$ (cf. e.{\,}g. \cite{RS}). The following result is now an immediate consequence of these observations.

\par

\begin{lemma}\label{basisM2}
Let $N$ be an integer. Then
$$
\{ (2|\alpha |+2d+1)^{-N}h_\alpha \} _{\alpha \in \mathbf N^d}
$$
is an orthonormal basis for $M^2_{2N}(\rr d)$.
\end{lemma}

\par

\subsection{Mapping properties of $\mathfrak V$ on
$M^2_N$}\label{subsec2.4}

\par

We shall now prove $A_{0,N}^2=A_N^2$ when $N$ is a non-zero even
integer. Important parts of these investigations are based upon the series
representation of analytic functions, using the fact that every $F\in A(\cc d)$
is equal to its Taylor series, i.{\,}e.
\begin{equation}\label{Taylorexp}
F(z)=\sum _{\alpha \in \mathbf N^d} a_\alpha \frac {z^\alpha}{(\alpha !)^{1/2}},\qquad
a_\alpha  = \frac {(\partial ^\alpha F)(0)}{(\alpha !)^{1/2}}.
\end{equation}
We also recall the result from \cite{B1} that $A_0^2(\cc d)=A^2(\cc d)$, and that $F\in A^2(\cc d)$, if and only if the coefficients in
\eqref{Taylorexp} satisfy
$$
\nm {(a_\alpha )_{\alpha \in \mathbf N^d}}{l^2}=\sum _{\alpha \in \mathbf N^d} |a_\alpha |^2 <\infty .
$$
Furthermore, $F=\mathfrak Vf\in A^2(\cc d)$ if and only if 
$f\in L^2(\rr d)$ satisfies
\begin{equation}\label{hermexp}
f(x)=\sum _{\alpha \in \mathbf N^d} a_\alpha h_\alpha (x),
\end{equation}
i.e., $f$ inherites the coefficients from $F$, and, since $\mathfrak{V}$ is isometric, 
\begin{equation}\label{normequiv}
\nm F{A^2} = \nm f{L^2}=\nm {(a_\alpha )_{\alpha \in \mathbf N^d}}{l^2}.
\end{equation}

\par

We now have the following result.

\par

\begin{prop}\label{imageM2N}
Let $N$ be an integer. Then the following is true:
\begin{enumerate}
\item $A_0(\sigma _{2N},L^2(\rr d))$ consists of all $F\in A(\cc d)$ with
expansion given by \eqref{Taylorexp}, where
\begin{equation}\label{normeq}
\nmm F\equiv \nm {(a_\alpha \eabs \alpha ^N)_{\alpha \in \mathbf N^d}}{l^2}<\infty .
\end{equation}
Furthermore, $\nmm \cdo $ and  $\nm \cdo {A(\sigma _{2N},L^2)}$  are
equivalent norms;

\vrum

\item $A(\sigma _{2N},L^2(\rr d)) = A_0(\sigma _{2N},L^2(\rr d))$.
\end{enumerate}
\end{prop}

\par

For the proof we recall that
\begin{equation}\label{BargHermmap}
(\mathfrak V h_\alpha )(z)=\frac {z^\alpha}{(\alpha !)^{1/2}}
\end{equation}
(cf. \cite {B1}) and that $\mathscr{S}_0(\rr d)$ is the set of all sums in \eqref{hermexp} such that $a_\alpha =0$ except for finite numbers of $\alpha$.

\par

\begin{proof}
(1) First we consider the case when $F\in P(\cc d)$, and we let
$a_\alpha$ be as in \eqref{Taylorexp}. Then it follows from \eqref{BargHermmap}
that $F$ is equal to $\mathfrak Vf$, where $f\in
\mathscr S_0(\rr d)$ is given by the finite sum \eqref{hermexp}. By \eqref{bargstft2}, 
Proposition \ref{identprop}, Lemma \ref{basisM2}, and \eqref{normequiv} it follows that
\begin{equation}\label{diskeqnormM2}
C^{-1}\nm F{A(\sigma _{2N},L^2)}\le \nm {((2|\alpha
|+2d+1)^Na_\alpha )_{\alpha}}{l^2}\le C\nm F{A(\sigma _{2N},L^2)},
\end{equation}
for some constant $C$ which is independent of $F\in P(\cc d)$. Since
$\mathscr S_0$ is dense in $M^2_{2N}$, it follows that
\eqref{diskeqnormM2} holds for each $F\in A_0(\sigma _{2N},L^2)$ when
$a_\alpha$ is given by \eqref{Taylorexp}. This proves (1).

\par

In order to prove (2) we recall that $\mathfrak{V}\, :\, M^2_{2N}\mapsto A_0(\sigma_{2N},L^2)$ is a bijective isometry in view of Proposition \ref{mainstep1prop}. Hence Lemma \ref{basisM2} 
together with \eqref{BargHermmap} show that
\begin{equation}\label{ONbasisA2N}
\left \{ (2|\alpha |+2d+1)^{-N}\frac {z^\alpha }{(\alpha !)^{1/2}}
\right \}
\end{equation}
is an orthonormal basis for $A_0(\sigma_{2N},L^2)$. By Proposition \ref{mainstep1prop} 
$\nm{F_0}{A_0(\sigma _{2N},L^2)}=\nm {F_0}{A(\sigma _{2N},L^2)}$ when $F_0\in A_0(\sigma _{2N},L^2)$.
Hence $A_0(\sigma _{2N},L^2)$ is a closed subspace of $A(\sigma _{2N},L^2)$. Consequently, we have the unique decomposition 
$$
A(\sigma _{2N},L^2) = A_0(\sigma _{2N},L^2)\oplus (A_0(\sigma _{2N},L^2))^\bot ,
$$
and it follows that \eqref{ONbasisA2N} is an orthonormal sequence in $A(\sigma _{2N},L^2)$. 
The fact that every $F\in A(\sigma _{2N},L^2)$ has a Taylor expansion now implies that \eqref{ONbasisA2N} is an orthonormal basis for $A(\sigma _{2N},L^2)$. Hence
$(A_0(\sigma _{2N},L^2))^\bot=\{0\}$ and the result follows.
\end{proof}

\par

\begin{cor}\label{imageSprime}
There is equality in \eqref{mapofSprime} in case $p=2$, i.e., 
$$
\sets {\mathfrak Bf}{f\in \mathscr S'(\rr d)} = \sets {F\in A(\cc d)}{\nm {Fe^{-|\cdo
|^2/2}\sigma _{-N}}{L^2}<\infty\ \text{for some}\ N\ge 0}.
$$
\end{cor}

\par

\begin{proof}
The result follows from Proposition \ref{imageM2N} and the fact that
$$
\cup _{N\in \mathbf Z}M^2_N(\rr d)=\mathscr S'(\rr d).
$$
\end{proof}

\par

\subsection{Mapping properties of $\mathfrak V$ on $\mathscr S'$, and
 proof of the main theorem}\label{subsec2.5}

\par

We shall now consider the relation \eqref{mapofSprime} and prove that
we indeed have equality when $1\le p\le 2$. In order to do this we
need the following lemma. Here we let $B_r(z)$ denote
the open ball in $\cc d$ with radius $r$ and center at $z\in \cc d$.

\par

\begin{lemma}\label{ballfamily}
There is a family $(B_j)_{j\in J}$ of open balls $B_j$ such that the
following conditions are fulfilled:
\begin{enumerate}
\item $\complement B_4(0)\subseteq \cup B_j$;

\vrum

\item $B_j=B_{r_j}(z_j)$ for some $r_j$ and $z_j$ such
that $|z_j|\ge 4$, $r_j\le 1/|z_j|$;

\vrum

\item there is a finite bound on the number of overlapping balls $B_{4r_j}(z_j)$.
\end{enumerate}
\end{lemma}

\par

\begin{proof}
Let $k\geq4$ and let $N$ be a large integer, and consider the spheres
$$
S_{k,l}=\{ z\in \cc d\, ;\ |z|=k+l/{kN}\},\quad l=0,\dots ,kN-1.
$$
On each sphere $S_{k,l}$, choose a finite number of points $z_j$
in such way that for any two closest points $z$ and $w$ the distance
between them is $1/2k\le |z-w|\le 1/(k+1)$. It is easily seen that such a
sequence $({z_j})$ exists when $N$ is chosen large
enough. The result now follows if we choose $B_j=B_{r_j}(z_j)$ with
$r_j=1/(k+1)$.
\end{proof}

\par

We have now the following result.

\par

\begin{prop}\label{imageSprime2}
Let $p\in [1,2]$ be fixed. Then there is equality in  
\eqref{mapofSprime}.
\end{prop}

\par

\begin{proof}
Let $\Omega _p$ be the set on the right-hand side of \eqref{mapofSprime}. In view of
Corollary \ref{imageSprime}, it suffices to prove that $\Omega _p$ is
independent of $p$. First assume that $p_1\le p_2$, and let $r\in
[1,\infty]$ be such that $1/p_2+1/r=1/p_1$. Then it follows from
H{\"o}lder's inequality that
\begin{multline*}
\nm {Fe^{-|\cdo |^2/2}\eabs
\cdo ^{-N-d-1}}{L^{p_1}} = \nm {(Fe^{-|\cdo |^2/2}\eabs
\cdo ^{-N})\eabs \cdo ^{-d-1}}{L^{p_1}}
\\[1ex]
\le C\nm {Fe^{-|\cdo |^2/2}\eabs
\cdo ^{-N}}{L^{p_2}},
\end{multline*}
where $C=\nm {\eabs \cdo ^{-d-1}}{L^r}<\infty$. This proves that
$\Omega _{p_2}\subseteq \Omega _{p_1}$.

\par

The result therefore follows if we prove that $\Omega _1\subseteq
\Omega _2$. Assume that $F\in \Omega _1$. It suffices to prove that
\begin{equation}\label{L2normest}
\int _{|z|\ge 4}|F(z)\eabs z^{-N}e^{-|z|^2/2}|^2\, d\lambda (z)<\infty ,
\end{equation}
for some $N\ge 0$. Here and in what follows, $d\lambda (z)$ denotes the Lebesgue measure
on $\cc d$.

\par

Since $F\in A(\cc d)$, the mean-value property for harmonic functions
gives
$$
F(z)=C|z|^{-d}\int _{|w|\le 1/|z|}F(z+w)\, d\lambda (w),
$$
where $1/C$ is the volume of the $d$-dimensional unit ball. Since
$$
C^{-1}e^{-|z|^2} \le e^{-|z+w|^2} \le Ce^{-|z|^2},\quad C^{-1}\eabs
z\le \eabs {z+w}\le C\eabs z\quad
\text{and}\quad \eabs z\le C|z|
$$
for some constant $C>0$, when $|w|\le 1/|z|$ and $|z|\ge 3$, we get
\begin{multline}\label{intest2}
\int _{|z|\ge 4}|F(z)\eabs z^{-N}e^{-|z|^2/2}|^2\, d\lambda (z)
\\[1ex]
\le C_1\int _{|z|\ge 4}\Big (\int _{|w|\le 1/|z|}|F(z+w)|\, d\lambda
(w)\eabs z^{-N-d}e^{-|z|^2/2}\Big )^2\, d\lambda (z)
\\[1ex]
\le C_2\int _{|z|\ge 4}\Big (\int _{|w|\le 1/|z|}|F(z+w) \eabs
{z+w}^{-N-d}e^{-|z+w|^2/2}|\, d\lambda (w)\Big )^2\, d\lambda (z)
\\[1ex]
=C_2 \int _{|z|\ge 4}\Big (\int _{|w-z|\le 1/|z|}|F(w) \eabs
{w}^{-N-d}e^{-|w|^2/2}|\, d\lambda (w)\Big )^2\, d\lambda (z).
\end{multline}

\par

Now let $B_j$ be as in Lemma \ref{ballfamily}. Then Lemma \ref{ballfamily} (1) gives that the integral on the right-hand side of \eqref{intest2} is estimated from above by
\begin{equation*}
C\sum _{j\in J} \int _{B_j}\Big ( \int _{|w-z|\le 1/|z|}|F(w)\eabs
w^{-N-d}e^{-|w|^2/2}|\, d\lambda (w) \Big )^2\eabs z^{2d}\, d\lambda
(z).
\end{equation*}
Since $|w-z_j|\le 4/|z_j|$ when $|w-z|\le 1/|z|$ and $z\in B_j$, the last integral can be estimated by
\begin{multline*}
C_1\sum _{j\in J} \int _{B_j}\Big ( \int _{|w-z_j|\le
4/|z_j|}|F(w)\eabs w^{-N-d}e^{-|w|^2/2}|\, d\lambda (w) \Big )^2\eabs
z^{2d}\, d\lambda (z)
\\[1ex]
\le C_2 \sum _{j\in J} \Big ( \int _{w\in B_{4r_j}(z_j)}|F(w)\eabs
w^{-N-d}e^{-|w|^2/2}|\, d\lambda (w) \Big )^2
\\[1ex]
\le C_2 \Big ( \sum _{j\in J} \int _{w\in B_{4r_j}(z_j)}|F(w)\eabs
w^{-N-d}e^{-|w|^2/2}|\, d\lambda (w) \Big )^2
\\[1ex]
\le C_3 \Big ( \int _{\cc d}|F(w)\eabs w^{-N-d}e^{-|w|^2/2}|\,
d\lambda (w) \Big )^2,
\end{multline*}
for some constants $C_1,\dots ,C_3$. Here the first inequality
follows from the fact that $\int _{B_j}\eabs z^{2d}\, d\lambda (z)\le
C$ for some constant $C$ which is independent of $j$ by the property (2) in Lemma \ref{ballfamily},  and the last two
inequalities follow from the fact that there is a finite number of of overlapping balls
$B_{4r_j}(z_j)$ by (3) in Lemma \ref{ballfamily}. Summing up we have proved that
$$
\Big ( \int _{|z|\ge 4}|F(z)\eabs z^{-N}e^{-|z|^2/2}|^2\, d\lambda (z)\Big )^{1/2}\le
C\nm {F\eabs \cdo ^{-N-d}e^{-|\cdo |/2}}{L^1},
$$
for some constant $C$. The proof is complete.
\end{proof}

\par

\begin{proof}[Proof of Theorem \ref{mainthm}]
By Proposition \ref{mainstep1prop} it follows that the map $f\mapsto
\mathfrak Vf$ is an isometric injection from $M(\omega ,\mathscr B)$
to $A(\omega ,\mathscr B)$. We have to show that this mapping is
surjective.

\par

Therefore assume that $F\in A(\omega ,\mathscr B)$. By Propositions
\ref{Binvemb}, \ref{imageM2N} and \ref{imageSprime2},
there is an element $f\in \mathscr S'(\rr d)$ such that $F=\mathfrak
Vf$. We have
$$
\nm f{M(\omega ,\mathscr B)} =\nm {\mathfrak Vf}{A(\omega ,\mathscr
B)} =\nm F{A(\omega ,\mathscr B)}<\infty .
$$
Hence, $f\in M(\omega ,\mathscr B)$, and the result follows. The proof
is complete.
\end{proof}

\par

\section{Some consequences}\label{sec3}

\par

In this section we present some results which are straight-forward
consequences of Theorem \ref{mainthm} and well-known properties for
modulation spaces. Most of these results can be found in
\cite{FGW,Gc1,Gc2,GW}.

\par

We start with introducing some notations. We set
$$
A^{p,q}_{(\omega )}(\cc d) = A(\omega ,L^{p,q}(\rr {2d}))\quad \text{and}\quad
A^{p}_{(\omega )} = A^{p,p}_{(\omega )},
$$
when $\omega \in \mathscr P(\cc d)$ and $p,q\in [1,\infty ]$. We also
set 
$$
A^{p,q} = A^{p,q}_{(\omega )}\quad \text{and}\quad A^{p} =
A^{p}_{(\omega )}\quad \text{when}\quad \omega =1.
$$

\par

Let
$$
d\mu (w) = \pi ^{-d}e^{-|w|^2}d\lambda (w),
$$
where $d\lambda (z)$ is the Lebesgue measure on $\cc d$. We recall
from \cite{B1,B2} that the standard scalar product on $A^{2}(\cc d)$ is
given by
\begin{equation}\label{Fochspaceform}
(F,G)_{A^{2}} \equiv \int _{\cc d} F(w)\overline
{G(w)}\, d\mu (w).
\end{equation}
Furthermore, there is a convenient reproducing kernel on $A_{\mathscr
S}'(\cc d)$, given by the formula
\begin{equation}\label{reproducing}
F(z)= \int _{\cc d} e^{(z,w)}F(w)\, d\mu (w),\quad F\in
A_{\mathscr S}'(\cc d),
\end{equation}
where $(\cdo ,\cdo )$ is the scalar product on $\cc d$
(cf. \cite{B1,B2}). For future references we observe that \eqref{reproducing} is the same as
\begin{equation}\tag*{(\ref{reproducing})$'$}
F(z) = \pi ^{-d}\scal {F\cdot e^{(z,\cdo )}}{e^{-|\cdo |^2}},\quad F\in
A_{\mathscr S}'(\cc d),
\end{equation}

\par

\subsection{Embedding and duality properties}

\par

We shall now discuss embedding properties. The following result follows
immediately from Proposition \ref{p1.4} (2) and Theorem \ref{mainthm}.

\par

\begin{prop}
Let $p_j,q_j\in [1,\infty ]$, $\omega
,\omega _j,v,v_0\in \mathscr P(\rr {2d})$ for $j=1,2$ be such that
$p_1\le p_2$, $q_1\le q_2$, $\omega$ is $v$-moderate, and $\omega
_2\le C\omega _1$ for some
constant $C$. Also let $\mathscr B$ be a translation invariant
BF-space on $\rr {2d}$ with respect to $v_0$. Then
\begin{alignat*}{2}
A_{\mathscr S}(\cc d) \subseteq A^{p_1,q_1}_{(\omega _1)}(\cc d)  
   &\subseteq A^{p_2,q_2}_{(\omega _2)}(\cc d)\subseteq
A_{\mathscr S}'(\cc d),
\\[1ex]
A^1_{(v_0v)}(\cc d)\subseteq A(&\omega ,\mathscr B)\subseteq
A^{\infty}_{(1/(v_0v))}(\cc d).
\end{alignat*}
\end{prop}

\par

\begin{prop}\label{dense}
  Let $\omega \in \mathscr P(\cc d)$. Then
  $P(\cc d)$ is dense in $A^{p,q}_{(\omega )}(\cc d)$ when $1\le
  p,q<\infty$.
\end{prop}

\par

\begin{proof}
Recall that $\mathscr S_0(\rr d)$, the set of finite linear
combinations of the Hermite functions, is dense in $\mathscr S(\rr d)$
(cf. \cite[Theorem V.13]{RS}). Hence the result follows immediately from Proposition \ref{p1.4},
Theorem \ref{mainthm}, and the fact that  $\mathfrak V(\mathscr
S_0(\rr d))=P(\cc d)$.
\end{proof}

\par

\begin{prop}
   Let $\omega \in \mathscr P(\rr {2d})$ and $p,q\in[1,\infty]$. Then
the form \eqref{Fochspaceform} on $P(\cc d)$ extends to a
   continuous sesquilinear form on $A^{p,q}_{(\omega )}(\cc d)\times
   A^{p',q'}_{(1/\omega )}(\cc d)$, and
   \begin{equation}\label{CSi}
      |(F,G)_{A^{2}}| \le \nm F{A^{p,q}_{(\omega )}} \nm
      G{A^{p',q'}_{(1/\omega )}}.
   \end{equation}
   This extension is unique, except when $p=q'\in \{ 1,\infty
\}$.

\par

Moreover, let
   \begin{equation}\label{eqn}
      \nmm F \equiv \sup |(F,G)_{A^{2}}| ,
   \end{equation}
    where the supremum is taken over all $G\in P(\cc d)$ (or $G\in
    A^{p',q'}_{(1/\omega )}(\cc d)$)
    such that $\nm G{A^{p',q'}_{(1/\omega )}}\le 1$. Then $\nmm \cdo$ and $\nm \cdo {A^{p,q}_{(\omega )}}$ are equivalent norms on $A^{p,q}_{(\omega
   )}(\cc d)$.
\end{prop}

\par

\begin{proof}
  The extension assertions and the inequality \eqref{CSi} are immediate
consequences of Proposition \ref{p1.4} (3), Theorem \ref{mainthm} and H{\"o}lder's inequality.

\par

  Let
\begin{align*}
\Omega _M &= \sets {g\in M^{p',q'}_{(1/\omega )}(\rr d)}{\nm
g{M^{p',q'}_{(1/\omega )}}\le 1},
\\[1ex]
\Omega _A &= \sets {G\in A^{p',q'}_{(1/\omega )}(\cc d)}{\nm
G{A^{p',q'}_{(1/\omega )}}\le 1}.
\end{align*}
For any $F\in A^{p,q}_{(\omega)}$ there is a unique $f\in
    M^{p,q}_{(\omega)}$ such that $\mathfrak Vf=F$. By Proposition
\ref{p1.4} (3) and Theorem \ref{mainthm} we get
   \begin{multline*}
     \|F\|_{A^{p,q}_{(\omega)}} = \|f\|_{M^{p,q}_{(\omega)}} 
      \le C \sup _{\substack{g\in \mathscr{S}\cap \Omega _M}} |(f,g)_{L^2}|
   \\[1ex] 
       \leq C\sup _{\substack{g\in \Omega _M}} |(f,g)_{L^2}|
       =C\sup _{\substack{G\in \Omega _A}}|(F,G)_{A^{2}}|
       \leq  C \|F\|_{A^{p,q}_{(\omega)}}, 
   \end{multline*}
   for some constant $C$, where the first  inequality follows from Proposition \ref{p1.4} 
   (3) and the last one from \eqref{CSi}.
   Since any $g\in\mathscr{S}$ can be approximated by its truncated
   Hermite expansion, the supremum over
   $\mathscr{S}$ may be substituted for a supremum over the finite
   Hermite expansions. 
   These, in turn, are the inverse images of the polynomials in $P(\cc
   d)$ which proves the last statement.
\end{proof}

\par

\begin{rem}
We note that the integral in \eqref{Fochspaceform} is well-defined
when $F\in A^{p,q}_{(\omega )}(\cc d)$, $G\in A^{p',q'}_{(1/\omega
)}(\cc d)$, $\omega \in \mathscr P(\cc d)$ and $p,q\in
[1,\infty]$. Also in the case $p=q'\in \{ 1,\infty \}$, we take this
integral as the definition of $(F,G)_{A^2}$, and we remark that the
extension of the form $(\cdo ,\cdo )_{A^2}$ on $P(\cc d)$ to
$A^{p,q}_{(\omega )}(\cc d)\times A^{p',q'}_{(1/\omega )}(\cc d)$ is
unique also in this case, if in addition narrow convergence is imposed
(cf Definition \ref{narrowdef} and Proposition \ref{narrowprop}
below).
\end{rem}

\par

\begin{prop}
   Let $\omega \in \mathscr P(\cc d)$ and $1\leq p,q<\infty$. Then the
   dual of $A^{p,q}_{(\omega)}(\cc d)$ can be identified with
   $A^{p',q'}_{(1/\omega)}(\cc d)$ through the form $(\cdo ,\cdo
   )_{A^{2}}$. Moreover, $P(\cc d)$ is weakly dense in $A^\infty
   _{(\omega )}(\cc d)$.
\end{prop}

\par

\begin{proof}
   The result is an immediate consequence of Proposition \ref{p1.4}
(4), Theorem \ref{mainthm}, and the fact that $\mathscr S_0$ is dense
in $\mathscr S$.

\end{proof}

\par

\subsection{Reproducing kernel and Bargmann-Toeplitz operators}

\par

For general $F\in L^2(d\mu )$, it is proved in \cite{B1} that the right-hand sides of \eqref{reproducing} and \eqref{reproducing}$'$ defines an orthonormal projection $\Pi _A$ of elements in $L^2(d\mu )$ onto $A^2(\cc d)$. We recall that $A^2$ is the image of $L^2$ under the Bargmann transform. In what follows we address equivalent projections where the Bargmann transform is replaced by the short-time Fourier transform. We use these relations to extend $\Pi _A$ to more general spaces of distributions.

\par

When dealing with the short time Fourier transform, it is convenient
to consider the twisted convolution $\widehat *$ on $L^1(\rr {2d})$,
which is defined by the formula
$$
(F\hatconv G)(x,\xi )=(2\pi )^{-d/2}\iint F(x-y,\xi -\eta
)G(y,\eta )e^{-i\scal {x-y}\eta}\, dyd\eta .
$$
(Cf. e.{\,}g. \cite{FG1,Gc2}.) By straight-forward computations it
follows that $\widehat *$
restricts to a continuous multiplication on $\mathscr S(\rr
{2d})$. Furthermore, the map $(F,G)\mapsto F\hatconv G$ from
$\mathscr S(\rr {2d})\times \mathscr S(\rr {2d})$ to $\mathscr S(\rr
{2d})$ extends uniquely to continuous mappings from $\mathscr S'(\rr
{2d})\times \mathscr S(\rr {2d})$ and $\mathscr S(\rr {2d})\times
\mathscr S'(\rr {2d})$ to $\mathscr S'(\rr {2d})\bigcap C^{\infty}(\rr
{2d})$.

\par

\begin{rem}\label{conseqtwisted}
By Fourier's inversion formula, it follows that
\begin{equation}\label{windowtransf}
(V_{\phi _1}f)\hatconv (V_{\phi _2}\phi _3) = (\phi _3,\phi
_1)_{L^2(\rr d)}\cdot V_{\phi _2}f
\end{equation}
for every $f\in \mathscr S'(\rr d)$ and every $\phi _j\in \mathscr
S(\rr d)$. The relation \eqref{windowtransf} is used in \cite{FG1,Gc2}
to prove the following properties:
\begin{enumerate}
\item The modulation spaces are independent of the choice of
window functions (cf. Proposition \ref{p1.4} (1));

\vrum

\item Let $\phi \in \mathscr S(\rr d)$ satisfy $\nm
{\phi}{L^2}=1$, and let $\Pi$ be the mapping on $\mathscr S'(\rr
{2d})$, given by
\begin{equation}\label{Pidef}
\Pi F \equiv F\hatconv (V_\phi \phi ).
\end{equation}
Also let $\mathscr B$ be a
translation invariant BF-space on $\rr {2d}$, $\omega \in \mathscr
P(\rr {2d})$, and set
$$
V_\phi (\Sigma )\equiv \sets {V_\phi f}{f\in 
\Sigma },
$$
when $\Sigma \subseteq \mathscr S'(\rr d)$. Then
\begin{alignat}{3}
&\Pi \, :\, & \mathscr S(\rr {2d}) &\to V_\phi (\mathscr S(\rr d)) &  & \subseteq \mathscr S(\rr {2d})\label{projections1}
\\[1ex]
&\Pi \, :\, & \mathscr S'(\rr {2d}) &\to V_\phi (\mathscr S'(\rr d)) &  & \subseteq \mathscr S' (\rr {2d})\label{projections2}
\\[1ex]
&\Pi \, :\, & \mathscr B(\omega ) &\to V_\phi (M(\omega ,\mathscr B)) & & \phantom =\label{projections3}
\end{alignat}
are continuous projections. Furthermore, if in
addition $\phi$ is given by \eqref{phidef}, then it follows by
straight-forward computations that $\Pi$ is self-adjoint on $L^2(\rr
{2d})$. Hence, for such a 
choice of $\phi$ it follows that $\Pi$ is an orthonormal
projection from  $L^2(\rr {2d})$ to $V_\phi (L^2(\rr d))$.
\end{enumerate}
\end{rem}

\par

Now we recall that the orthonormal $\Pi _A$ of $L^2(d\mu )$ onto $A^2(\cc d)$ is given by the right-hand sides of the reproducing formulas \eqref{reproducing} and \eqref{reproducing}$'$, i.{\,}e. 
\begin{equation}\label{A2projection}
(\Pi _AF)(z)= \int _{\cc d} e^{(z,w)}F(w)\, d\mu (w),\quad F\in
L^2(d\mu ).
\end{equation}
We extend the definition of $\Pi _A$ to the set
$$
\mathfrak S'(\cc d) \equiv \sets {F\in \mathscr D'(\cc d)}{Fe^{-|\cdo
|^2/2}\in \mathscr S'(\cc d)},
$$
by the formula
\begin{equation}\tag*{(\ref{A2projection})$'$}
(\Pi _AF)(z)= \scal {F\cdot e^{(z,\cdo )}}{e^{-|\cdo |^2}},\quad F\in
\mathfrak S'(\cc d),
\end{equation}
and we note that \eqref{A2projection}$'$ agree with \eqref{A2projection} when $F\in L^2(d\mu )$.

\par

We note that the set $\mathfrak S'(\cc d)$ is equal to $U_{\mathfrak V}(\mathscr S'(\rr {2d}))$, where $U_{\mathfrak V}$ is given by \eqref{UVdef}. Furthermore, by letting $\phi _j(x)=\phi (x)=\pi ^{-d/4}e^{-|x|^2/2}$, the reproducing formulas \eqref{reproducing} and \eqref{reproducing}$'$ are
straight-forward consequence of \eqref{bargstft2} and \eqref{windowtransf}. From these computations it also follows that $\Pi _A$ is the conjugation of $\Pi$ in \eqref{Pidef} by $U_{\mathfrak V}$, i.{\,}e.
\begin{equation}\label{conjugation}
\Pi _A = U_{\mathfrak V}\circ \Pi \circ U_{\mathfrak V}^{-1}.
\end{equation}

\par

The following result is now an immediate consequence of
these observations, Theorem \ref{mainthm} and \eqref{projections1}--\eqref{projections3}. Here we let
$$
\mathfrak S(\cc d) \equiv \sets {F\in \mathscr D'(\cc d)}{Fe^{-|\cdo
|^2/2}\in \mathscr S(\cc d)},
$$
which is the same as $U_{\mathfrak V}(\mathscr S(\rr {2d}))$.

\par

\begin{prop}\label{projectionprop}
Let $\mathscr B$ be a translation invariant BF-space on $\rr {2d}$,
and let $\omega \in \mathscr P(\rr {2d})$. Then the following hold:
\begin{enumerate}
\item $\Pi _A$ is a continuous projection from $\mathfrak S'(\cc d)$ to
$A_{\mathscr S}'(\cc d)$;

\vrum

\item $\Pi _A$ restricts to a continuous projection from $\mathscr
B_{\mathfrak V}(\omega )$ to $A(\omega ,\mathscr B)$;

\vrum

\item $\Pi _A$ restricts to a continuous projection from $\mathfrak S(\cc d)$
to $A_{\mathscr S}(\cc d)$.
\end{enumerate}
\end{prop}

\medspace

Next we consider Toeplitz operators in the context of Bargmann transform. 
It follows from \eqref{toeplitzdef} that if $a\in \mathscr S'(\rr {2d})$ and $f,\phi \in \mathscr S(\rr d)$, then
\begin{equation}\label{toeplitz2}
(V_\phi \circ \tp _\phi (a))f = \Pi (a\cdot F_0),\quad \text{where}\quad F_0=V_\phi f .
\end{equation}

\par

The close relation between the short-time Fourier transform and the Bargmann transform motivates the following definition.

\par

\begin{defn}
Let $a\in \mathscr S'(\cc d)$, and let $S$ be as in \eqref{Sdef}. Then
the Bargmann-Toeplitz operator $\operatorname T_{\mathfrak V}(a )$ is
the continuous operator on $A_{\mathscr S}'(\cc d)$, given by the
formula
$$
\operatorname T_{\mathfrak V}(a )F = \Pi _A((S^{-1}a )F).
$$
\end{defn}

\par

It follows from \eqref{conjugation} and \eqref{toeplitz2} that
$$
\operatorname T_{\mathfrak V}(a ) \circ \mathfrak V = \mathfrak V \circ \tp (a).
$$
The following result is now an immediate consequence of the latter property
and \cite[Theorem 3.1]{GT}. We recall that $\mathscr P_0$
consists of all smooth elements $\omega$ in $\mathscr P$ such that
$(\partial ^\alpha \omega )/\omega \in L^\infty$.

\par

\begin{prop}
Let $\omega \in \mathscr P(\cc d)$, $\omega _0\in \mathscr P_0(\cc
d)$, and let $\mathscr B$ be a translation invariant BF-space. Then
$\operatorname T_{\mathfrak V}(\omega _0)$ is continuous and bijective
from $A(\omega ,\mathscr B)$ to $A(\omega /\omega _0,\mathscr B)$.
\end{prop}

\par

\subsection{Mapping properties of Harmonic oscillator on modulation spaces}

\par

We shall now show how our investigations can be used to get spectral properties of harmonic oscillator. For each $t\in \mathbf C$, we let the ($t$-)harmonic oscillator be defined by
\begin{equation}\label{harmosc2}
H_t\equiv |x|^2 -\Delta +t,
\end{equation}
and we observe that this definition agree with \eqref{harmosc} when $t=d+1$. By \cite{B1} it follows that
\begin{equation}\label{Httransf}
\mathfrak V(H_tf)=2\left (\sum _{j=1}^d z_j\frac {\partial F}{\partial z_j}\right ) +d\cdot F,\quad F\in A_{\mathscr S}'(\cc d),
\end{equation}
which implies that if $F\in A_{\mathscr S}'(\cc d)$ is given by
\begin{equation}\label{Fserietransfor}
F(z) = (\mathfrak Vf)(z) = \sum _\alpha a_\alpha \frac {z^\alpha}{(\alpha !)^{1/2}},
\end{equation}
then
$$
\mathfrak V(H_t^Nf)(z) = \sum _\alpha a_\alpha \frac {(2|\alpha |+d+t)^Nz^\alpha}{(\alpha !)^{1/2}},
$$
as $N\ge 0$ is an integer.

\par

For the harmonic oscillator we have now the following result.

\par

\begin{thm}\label{thmharmonic}
Let $t\in \mathbf C$ and $N\in \mathbf N$ be fixed, $\omega \in \mathscr P(\rr d)$, and let $\mathscr B$ be a translation invariant BF-space on $\rr {2d}$. Then the following conditions are equivalent:
\begin{enumerate}
\item $H_t^N$ is continuous and bijective on $\mathscr S(\rr d)$;

\vrum

\item $H_t^N$ is continuous and bijective on $\mathscr S'(\rr d)$;

\vrum

\item $H_t^N$ is continuous and bijective from $L^2(\rr d)$ to $M^2_{-2N}(\rr d)$;

\vrum

\item $H_t^N$ is continuous and bijective from $M(\omega ,\mathscr B)$ to $M(\omega /\sigma _{2N},\mathscr B)$;

\vrum

\item $t\notin \sets {-d-2n}{n\in \mathbf N}$.
\end{enumerate}

\par

\noindent
Furthermore, if {\rm{(5)}} is fulfilled, then {\rm{(1)}}--{\rm{(4)}} hold for each $N\in \mathbf Z$.
\end{thm}

\par

By \eqref{Httransf} and Theorem \ref{mainthm} it follows that Theorem \ref{thmharmonic} is equivalent to the following proposition.

\par

\begin{prop}\label{propharmonic}
Let $t\in \mathbf C$ and $N\in \mathbf N$ be fixed, $\omega \in \mathscr P(\rr d)$, and let $\mathscr B$ be a translation invariant BF-space on $\rr {2d}$. Also let
$$
T=\left (2\left (\sum _{j=1}^d z_j\frac {\partial F}{\partial z_j}\right ) +d \right ) ^N
$$
Then the following conditions are equivalent:
\begin{enumerate}
\item $T$ is continuous and bijective on $A_{\mathscr S}(\cc d)$;

\vrum

\item $T$ is continuous and bijective on $A_{\mathscr S}'(\cc d)$;

\vrum

\item $T$ is continuous and bijective from $A^2(\cc d)$ to $A^2_{-2N}(\cc d)$;

\vrum

\item $T$ is continuous and bijective from $A(\omega ,\mathscr B)$ to $A(\omega /\sigma _{2N},\mathscr B)$;

\vrum

\item $t\notin \sets {-d-2n}{n\in \mathbf N}$.
\end{enumerate}
\end{prop}

\par

\begin{proof}[Proof of Theorem \ref{thmharmonic} and Proposition \ref{propharmonic}]
The results follow if we prove that each one of (1)--(4) implies (5) in Proposition \ref{propharmonic}, that (5) implies (3) in Proposition \ref{propharmonic}, and that (3) $\Rightarrow$ (4), (4) $\Rightarrow$ (1) and (4) $\Rightarrow$ (2) in Theorem \ref{thmharmonic}.

\par

First we assume that (5) in Proposition ref{propharmonic} does not hold. Then $Tz^\alpha=0$ when $|\alpha |=-(t+d)/2\in \mathbf N$. This implies that (1)--(4) in Proposition \ref{propharmonic} do not hold. It remains to prove that
(5) $\Rightarrow$ (3) in Proposition \ref{propharmonic}, and that (3) $\Rightarrow$ (4), (4) $\Rightarrow$ (1) and (4) $\Rightarrow$ (2) in Theorem \ref{thmharmonic}.

\par

Therefore, assume that (5) in Proposition \ref{propharmonic} holds, and let $S$ be the operator on $A_{\mathscr S}'(\cc d)$ which maps $F$ in \eqref{Fserietransfor} into
$$
\sum _\alpha a_\alpha \frac {(2|\alpha |+2d+1)^N}{(2|\alpha |+d+t)^N)}\frac {z^\alpha}{(\alpha !)^{1/2}}.
$$
Then $S\circ H_t=H$, where $H$ is given by \eqref{harmosc}. Furthermore, it follows from the assumptions on $t$ that $S$ is continuous and bijective on each $A^2_{2N_0}$ for every integer $N_0$. The assertion (3) in Proposition \ref{propharmonic} and Theorem \ref{thmharmonic} are now consequences of Proposition \ref{identprop}, and (4) in Theorem \ref{thmharmonic} follows from \cite[Theorem 3.1]{GT} and Theorem \ref{mainthm}.

\par

Finally, (1) and (2) in Theorem \ref{thmharmonic} now follows from (4) and the relations
$$
\mathscr S'(\rr d)=\cup _{\omega \in \mathscr P}M(\omega ,\mathscr B)\quad \text{and}\quad \mathscr S(\rr d)=\cap _{\omega \in \mathscr P}M(\omega ,\mathscr B).
$$
The proof is complete.
\end{proof}

\par

\begin{rem}
Let $N\ge 0$ be an integer, $t\in \mathbf C\setminus \sets {-d-2n}{n\in \mathbf N}$, $m,s\in \mathbf R$, and let $\operatorname{Sh}^{m}_{1}(\rr {2d})$ be the Shubin-class of all smooth symbols $a$ on $\rr {2d}$ which satisfy
$$
|\partial _x^\beta \partial _\xi^\alpha a(x,\xi )|\le C_{\alpha ,\beta}\eabs {x,\xi }^{m-|\alpha |-|\beta |},
$$
for some constants $C_{\alpha ,\beta}$ which only depend on $a$, $\alpha$ and $\beta$. Also let the pseudo-differential operator $\op _s(b)$ with symbol $b\in \mathscr S'(\rr {2d})$ be defined in the usual way (cf. e.{\,}g. \cite{GT}). Then it follows that $H_t^N=\op _s(a)$ for some $a\in \operatorname{Sh} ^{2N}_{1}(\rr {2d})$. Furthermore, since $\operatorname{Sh}^{0}_{1}(\rr {2d})$ is a Wiener algebra (cf. \cite{Beals,Bo1}), the proof of \cite[Theorem 2.1]{GT} in combination with Theorem \ref{thmharmonic} show that the inverse $H_t^{-N}$ of $H_t^N$ is equal to $\op _s(b)$ for some $b\in \operatorname {Sh}^{-2N}_1(\rr {2d})$.
\end{rem}

\par

\subsection{The narrow convergence}

\par

We shall now discuss the narrow convergence for modulation spaces and discuss corresponding concept in context of generalized Bargmann-Foch spaces.

\par

The main reason why introducing the narrow convergence in context of Bargmann-Foch spaces is to improve the
possibilities for approximating elements in $A^{p,q}_{(\omega )}(\cc d)$ with elements in $P(\cc d)$. In terms of norm convergence, Proposition \ref{dense} does not
guarantee that such approximations are possible when $p=\infty$ or $q=\infty$. In the
case $p=q'\notin \{ 1,\infty \}$, the
situation is usually handled by using weak$^*$-topology, if
necessary. However, the remaining case $p=q'\in \{ 1,\infty \}$
may be critical since $P(\cc d)$ is neither dense in $A^{\infty ,1}_{(\omega )}(\cc d)$ nor in 
$A^{1,\infty }_{(\omega )}(\cc d)$. Here we shall see that such
problems may be avoided by using narrow
convergence from \cite{To5} for modulation spaces in the $A^{p,q}_{(\omega )}$ spaces.

\par

First we define the narrow convergence of such spaces.

\par

\begin{defn}\label{narrowdef}
Let $\omega \in \mathscr P(\cc d)$, $S$ be as in \eqref{Sdef}, $p,q\in [1,\infty ]$ and let
$F_j,F\in A_{(\omega )}^{p,q}(\cc d)$, $j\ge 1$. Then $F_j$ is
said to converge narrowly to $F$ as $j$ turns to infinity whenever
\begin{enumerate}
\item $F_j\to F$ in $A_{\mathscr S}'(\cc d)$ as $j\to \infty$;

\vrum

\item if
\begin{align*}
H_j(\xi ) &= \Big (\int _{\rr d} |F_j(z)e^{-|z|^2/2}(S^{-1}\omega )(z)|^p\, dx\Big )^{1/p},\quad
\\[1ex]
H(\xi ) &= \Big (\int _{\rr d}|F(z)e^{-|z|^2/2} (S^{-1}\omega )(z)|^p\, dx\Big )^{1/p},
\end{align*}
with $z=x+i\xi$ and $x,\xi \in \rr d$, then $H_j\to H$ in $L^q(\rr d)$.
\end{enumerate}
\end{defn}

\par

The following proposition justifies the definition of narrow convergence.

\par

\begin{prop}\label{narrowprop}
Let $\omega \in \mathscr P(\cc d)$ and let $G\in
A^{1,\infty}_{(1/\omega )}(\cc d)$. Then the following hold:
\begin{enumerate}
\item $P(\cc d)$ is dense in $A^{\infty ,1}_{(\omega )}(\cc d)$ with
respect to narrow convergence;

\vrum

\item if $F_j\in A^{\infty ,1}_{(\omega )}(\cc d)$ converges narrowly
to $F\in A^{\infty ,1}_{(\omega )}(\cc d)$ as $j\to \infty$, then
$(F_j,G)\to (F,G)$ as $j\to \infty$.
\end{enumerate}
\end{prop}

\par

\begin{proof}
The result follows immediately from Proposition 1.10 and Lemma 1.11 in
\cite{To5}, Theorem \ref{mainthm} and the fact that $P(\cc d)$ is dense in
$A_{\mathscr S}(\cc d)$.
\end{proof}

\medspace


\begin{thebibliography}{150}
\bibitem{B1} V. Bargmann, \emph{On a Hilbert space of analytic
functions and an associated integral transform}, Comm. Pure
Appl. Math., \textbf{14} (1961), 187--214.

\bibitem{B2} \bysame \emph{On a Hilbert space of analytic functions
 and an associated integral transform. Part II. A family of related
function spaces. Application to distribution theory.}, Comm. Pure
Appl. Math., \textbf{20} (1967), 1--101.

\bibitem{Beals} R. Beals \emph{Characterization of pseudodifferential operators and applications},
Duke Math. J. \textbf{44} (1977), 45--57. 

\bibitem{Berezin71}
F.~A. Berezin.
\newblock Wick and anti-{W}ick symbols of operators.
\newblock {\em Mat. Sb. (N.S.)}, 86(128):578--610, 1971.

\bibitem{Bo1} J.-M. Bony \emph{Caract{\'e}risations des op{\'e}rateurs pseudo-diff{\'e}rentiels \rm{in: S{\'e}minaire sur les {\'E}quations aux D{\'e}riv{\'e}es Partielles, 1996--1997}}, Exp. No. XXIII., {\'E}cole Polytech., Palaiseau, 1997.

\bibitem{Dau} {I. Daubechies} \emph{Time-frequency localization
  operators: a geometric phase space approach}, IEEE
  Trans. Inform. Th. (4) \textbf{34}, (1988), 605--612.

\bibitem{F1}  H.~G.~Feichtinger \emph{Modulation spaces on locally
compact abelian groups. Technical report}, {University of
Vienna}, Vienna, 1983; also in: M. Krishna, R. Radha,
S. Thangavelu (Eds) Wavelets and their applications, Allied
Publishers Private Limited, NewDehli Mumbai Kolkata Chennai Hagpur
Ahmedabad Bangalore Hyderbad Lucknow, 2003, pp. 99--140.

\bibitem{F2} \bysame \emph{Wiener amalgams over Euclidean spaces and
some of their applications},
in: Function spaces (Edwardsville, IL, 1990), Lect. Notes in pure and
appl. math., 136, Marcel Dekker, New York, 1992, pp. 123--137.

\bibitem{Feichtinger5} \bysame \emph{Gabor frames and time-frequency
analysis of distributions}, {J. Functional
Anal. (2)} \textbf {146} (1997), 464--495.

\bibitem{F3} \bysame \emph{Modulation spaces: Looking back and ahead},
Sampl. Theory Signal Image Process. \textbf{5} (2006), 109--140.

\bibitem{FG1}  {H. G. Feichtinger and K. H. Gr{\"o}chenig}
\emph{Banach spaces related to integrable group representations and
their atomic decompositions, I}, J. Funct. Anal., \textbf{86}
(1989), 307--340.

\bibitem{FG2} \bysame \emph{Banach spaces related to
integrable group representations and their atomic decompositions, II},
Monatsh. Math., \textbf{108} (1989), 129--148.

\bibitem{FGW}  {H. G. Feichtinger and K. H. Gr{\"o}chenig, D. Walnut}
\emph{Wilson bases and modulation spaces}, {Math. Nach.} \textbf{155}
(1992), 7--17.

\bibitem{Fo}  {G. B. Folland} \emph
{Harmonic analysis in phase space}, {Princeton U. P., Princeton},
1989.

\bibitem{Grobner}  {P. Gr{\"o}bner} {Banachr{\"a}ume Glatter
Funktionen und Zerlegungsmethoden, Thesis}, \textit{University
of Vienna}, Vienna, 1992.

\bibitem{Gc1} {K. H. Gr{\"o}chenig} \emph {Describing
functions: atomic decompositions versus frames},
{Monatsh. Math.},\textbf{112} (1991), 1--42.

\bibitem{Gc2} K. Gr\"{o}chenig, \newblock \textit{Foundations of
Time-Frequency Analysis},
\newblock Birkh{\"a}user, Boston, 2001.

\bibitem{GL} K. Gr{\"o}chenig, M. Leinert \emph{Wiener's lemma
for twisted convolution and Gabor frames}, J. Amer. Math. Soc.,
\textbf{17} (2004), 1--18.

\bibitem{GT} K. Gr{\"o}chenig, J. Toft \emph{Isomorphism properties of
Toeplitz operators and pseudo-differential operators between
modulation spaces}, J. Math. An. (to appear), also available at arXiv, arXiv:0905.4954.

\bibitem{GW} K. Gr{\"o}chenig, D. Walnut \emph{A Riesz basis for
Bargmann-Fock space related to sampling and interpolation},
Ark. Mat. \textbf{30} (1992), 283--295.

\bibitem{He1} F. H{\' e}rau \emph{Melin--H{\"o}rmander inequality in a
Wiener type pseudo-differential algebra}, Ark. Mat., \textbf{39}
(2001), 311--38.

\bibitem{HTW} A. Holst, J. Toft, P. Wahlberg \emph{Weyl product
algebras and modulation spaces}, {J. Funct. Anal.}, \textbf{251}
(2007), 463--491.

\bibitem{Ho1}  L. H{\"o}rmander \emph{The Analysis of Linear
Partial Differential Operators}, vol {I--III},
Springer-Verlag, Berlin Heidelberg NewYork Tokyo, 1983, 1985.

\bibitem{JPR} S. Janson, J. Petree, R. Rochberg  \emph{Hankel forms
and the Fock space}, Rev. Mat. Iberoamericana \textbf{3} (1987),
61--138.

\bibitem{Le} {N. Lerner} \emph{The Wick calculus of pseudo-differential operators and some of its applications}, CUBO, \textbf {5} (2003), 213--236. 

\bibitem{RS}  {M. Reed, B. Simon} \emph{Methods of modern
mathematical physics}, Academic Press, London New York, 1979.

\bibitem{Sugimoto1} {M. Sugimoto, N. Tomita} \emph{The dilation
property of modulation spaces and their inclusion relation with Besov
Spaces}, {J. Funct. Anal. (1)}, \textbf{248} (2007),
79--106.

\bibitem{Toft2} J. Toft \emph{Continuity properties for
modulation spaces with applications to pseudo-differential calculus,
I}, {J. Funct. Anal. (2)}, \textbf{207} (2004),
399--429.

\bibitem{To04B} \bysame \emph{Convolution and embeddings for
weighted modulation spaces {\textit {in: P. Boggiatto, R. Ashino,
M. W. Wong (Eds)}} Advances in Pseudo-Differential Operators,}
Operator Theory: Advances and Applications \textbf{155},
Birkh{\"a}user Verlag, Basel 2004, pp. 165--186.


\bibitem{To5} \bysame \emph{Continuity
properties for modulation spaces with applications to
pseudo-differential calculus, II}, {Ann. Global Anal. Geom.},
\textbf{26} (2004), 73--106.

\bibitem{To6} \bysame \emph{Minizization under entropy conditions, with applications in lower bound problems}, {J. Math. Phys.}, \textbf{45} (2004), 3216--3227.

\bibitem{To8} \bysame \emph{Continuity and Schatten
properties for pseudo-differential operators on modulation spaces {\rm
{in: J. Toft, M. W. Wong, H. Zhu (eds)}} Modern Trends in
Pseudo-Differential Operators,} Operator Theory: Advances and
Applications, Birkh{\"a}user Verlag, Basel, 2007, 173--206.
\end{thebibliography}
\end{document}